\documentclass{amsart}
\usepackage{amssymb}
\newtheorem{theorem}{Theorem}[section]
\newtheorem{lemma}[theorem]{Lemma}

\newtheorem{fact}[theorem]{Fact}

\theoremstyle{definition}

\theoremstyle{remark}
\newtheorem{remark}[theorem]{Remark}

\numberwithin{equation}{section}

\renewcommand{\dim}{\mathrm{dim}}

\newcommand{\dens}{\mathrm{dens}}

\renewcommand{\span}{\mathrm{span}}
 
\newcommand{\dist}{\mathrm{dist}}
\newcommand{\conv}{\mathrm{conv}}

\newcommand{\R}{\mathbb{R}}
\newcommand{\N}{\mathbb{N}}

\newcommand{\X}{\mathrm{X}}
\renewcommand{\H}{\mathrm{H}}
\newcommand{\Y}{\mathrm{Y}}
\newcommand{\Z}{\mathrm{Z}}
\newcommand{\B}{\mathbf{B}}
\newcommand{\I}{\mathbf{I}}
\newcommand{\sign}{\mathrm{sign}}

\renewcommand{\S}{\mathbf{S}}

\renewcommand{\mod}{/}

\renewcommand{\ae}{\stackrel{m}{\sim}}

\begin{document}

\title{Convex-transitivity and function spaces}

%    Information for first author
\author{Jarno Talponen}
%    Address of record for the research reported here
\address{University of Helsinki, Department of Mathematics and Statistics, Box 68, (Gustaf H\"{a}llstr\"{o}minkatu 2b) FI-00014 University 
of Helsinki, Finland}
\email{talponen@cc.helsinki.fi}
%    \thanks will become a 1st page footnote.

%    General info
\date{\today}

\begin{abstract}
It is shown that if $\X$ is a convex-transitive Banach space and\\ 
$1\leq p\leq \infty$, then the closed linear span of the 
simple functions in the Bochner space $L^{p}([0,1],\X)$ is also convex-transitive. If $\H$ is an infinite-dimensional Hilbert space 
and $C_{0}(L)$ is convex-transitive, then $C_{0}(L,\H)$ is convex-transitive.
Some new fairly concrete examples of convex-transitive spaces are provided.
\end{abstract}

\maketitle

\section{Introduction}
This paper is concerned with several aspects of symmetries of real Banach spaces. We denote the closed unit ball of 
a real Banach space $\X$ by $\B_{\X}$ and the unit sphere of $\X$ by $\S_{\X}$.
A Banach space $\X$ is called \emph{transitive} if for each $x\in \S_{\X}$ the orbit
$\mathcal{G}_{\X}(x)\stackrel{\cdot}{=}\{T(x)|\ T\colon \X\rightarrow \X\ \mathrm{is\ an\ isometric\ automorphism}\}=\S_{\X}$. 
If $\overline{\mathcal{G}_{\X}(x)}=\S_{\X}$ (resp. $\overline{\conv}(\mathcal{G}_{\X}(x))=\B_{\X}$) 
for all $x\in\S_{\X}$, then $\X$ is called \emph{almost transitive} (resp. \emph{convex-transitive}). 

These concepts are motivated by the \emph{Banach-Mazur rotation problem} appearing in \cite[p.242]{Ba} which remains unsolved
to this day and asks whether each separable transitive space is isometrically a Hilbert space. A branch of the isometric theory 
of Banach spaces has developed around this problem, and we refer to \cite{BR2} for an extensive survey of this field. In this article we 
will investigate the abovementioned transitivity conditions, mainly convex-transitivity, in different settings.
For some important results related to this investigation, see \cite{GJK}, \cite{GR}, \cite{Kawamura} and \cite{Rambla}.

The main line of study in this paper involves the rotations of certain Banach-valued function spaces. 
It is already a classical result that the space $L^{p}$ is almost transitive for $p\in [1,\infty)$ and convex-transitive for 
$p=\infty$ (see \cite{Rol}). In \cite{GJK} the almost transitivity of the Bochner space $L^{p}([0,1],\X)$ 
was established in the case where $\X$ is almost transitive and $1\leq p<\infty$. 
See \cite{AR} for the analogous study of the spaces $C_{0}(L,\X)$. We will extend these investigations into the convex-transitive setting. 
We will show that if $\X$ is convex-transitive, then $L_{s}^{\infty}(\X)$ is convex-transitive. Here $L_{s}^{\infty}(\X)$ denotes 
the closed linear span of the simple functions in $L^{\infty}([0,1],\X)$. It also turns out that if $L$ is a locally 
compact space such that $C_{0}(L)$ is convex-transitive and $\H$ is an infinite-dimensional Hilbert space, then 
$C_{0}(L,\H)$ is convex-transitive. The arguments are different compared to \cite{GJK} and \cite{AR}, since, for example, the rotation
group $\mathcal{G}_{L^{\infty}}$ differs considerably from $\mathcal{G}_{L^{p}}$ for $p<\infty$. Some implications
of the transitivity conditions of $C_{0}(L,\X)$ for the topology of $L$ are also discussed. 

We will also give some new, fairly concrete examples of almost transitive and convex-transitive spaces.
These examples are built starting from classical function spaces. For example, it turns out that 
the subspace of $L^{1}$ consisting of all the functions whose average is $0$, is almost transitive in the usual norm. 

Finally, we will study non-separable transitive spaces. Recall that any Banach space $\X$ of density character $\kappa$ 
can be isometrically embedded as a subspace of an almost transitive space of the same density character 
(see \cite{Lu} and \cite[p.8-11]{BR2}) or of a transitive space (see \cite[Cor 2.21]{BR2}). 
Starting from a transitive space one can also pass to a separable almost transitive subspace 
(see \cite{Ca4} and also \cite[Thm. 2.24]{BR2}). We will prove here that each transitive 
space contains a transitive subspace, the density character of which does not exceed $2^{\aleph_{0}}$. 

\subsection*{Preliminaries}
Throughout this article we will consider \emph{real} Banach spaces denoted by $\X,\Y$ and $\Z$ unless otherwise stated.
For non-empty subsets $A,B\subset \X$ we put $\mathrm{dist}(A,B)=\inf_{a\in A, b\in B}||a-b||$. The group of rotations $\mathcal{G}_{\X}$
of $\X$ consists of isometric automorphisms $T\colon\X\rightarrow\X$, the group operation being the composition
of the maps and the neutral element is the identity map $\I\colon \X\rightarrow\X$. If $f\in \X^{\ast}$ and $x\in \X$, then
we denote by $f\otimes x$ the map $f(\cdot)x\colon \X\rightarrow [x]$. Recall that $x\in \S_{\X}$ is called a \emph{big point} if
$\overline{\conv}(\mathcal{G}_{\X}(x))=\B_{\X}$.

We refer to \cite{Lac} for background on measure algebras and 
isometries of $L^{p}$-spaces. In what follows $\Sigma$ is the completed sigma algebra of Lebesgue measurable sets on $[0,1]$. We denote by 
$m\colon \Sigma\rightarrow \R$ the Lebesgue measure. Let us define an equivalence relation $\ae$ on $\Sigma$
by setting $A\ae B$ if $m((A\cup B)\setminus (A\cap B))=0$. 

For an introduction to ordinal numbers and such matters we refer to \cite{En}. 
We denote the cardinality of a set $A$ by $|A|$. Here we will apply cardinal arithmetic notations, so that $2^{\aleph_{0}}=|\R|$. 
Recall that the density character of $\X$, $\dens(\X)$ for short, is the smallest cardinal $\kappa$ such that 
$\X$ contains a dense subset of cardinality $\kappa$. Given a limit ordinal $\kappa$, its cofinality is 
$\mathrm{cf}(\kappa)\stackrel{\cdot}{=}\min\{|A|:\ A\subset\kappa,\ \sup A=\kappa\}$.

\section{Convex-transitivity of Banach-valued function spaces}

For notational simplicity we abbreviate the Lebesgue-Bochner space 
$L^{p}([0,1],\X)$ by $L^{p}(\X)$ for $1\leq p\leq\infty$. 
Recall that this space consists of strongly measurable maps $f\colon [0,1]\rightarrow \X$ endowed with the norm 
\[||f||_{L^{p}(\X)}^{p}=\int_{0}^{1}||f(t)||_{\X}^{p}\ \mathrm{d}t,\quad \mathrm{for}\ p\in [1,\infty)\]
and $||f||_{L^{\infty}(\X)}=\underset{t\in [0,1]}{\mathrm{ess\ sup}}||f(t)||_{X}$. 
It is not hard to verify that the subspace 
\[L_{s}^{\infty}(\X)=\overline{\span}(\{\chi_{A}x|A\in\Sigma,\ x\in \X\})\subset L^{\infty}(\X)\] 
is non-separable if $\dim(\X)>0$ and $L_{s}^{\infty}(\X)\subsetneq L^{\infty}(\X)$ if $\dim(\X)=\infty$. 
We refer to \cite{DU} for precise definitions and background information regarding the Banach-valued function spaces
appearing in this section.  

\subsection{Convex-transitivity of $L^{p}(\X)$}
Recall that $L^{\infty}$ is convex-transitive (see \cite{PR} and \cite{Rol}). 
Greim, Jamison and Kaminska proved that $L^{p}(\X)$ is almost transitive if $\X$ is almost transitive and $1\leq p<\infty$, 
see \cite[Thm. 2.1]{GJK}. We will present an analogous result for convex-transitive spaces, that is, 
if $\X$ is convex-transitive, then $L_{s}^{\infty}(\X)$ and $L^{p}(\X)$ are also convex-transitive for $1\leq p<\infty$. 

Note that $L^{p}(\ell^{p})=L^{p}$ isometrically, where $\ell^{p}$ is not convex-transitive for $p\neq 2$, 
so that the above implications cannot be reversed. 
We know neither an example of a convex-transitive space $\X$ such that $L^{\infty}(\X)$ fails convex-transitivity, nor an example of a 
transitive space $\Y$ such that $L^{2}(\Y)$ fails transitivity. 

\begin{theorem}\label{th: LPX} 
Let $\X$ be a convex-transitive space. Then $L_{s}^{\infty}(\X)$ and $L^{p}(\X)$, where $1\leq p<\infty$, are convex-transitive.
\end{theorem}
We will make some preparations before giving the proof. The following two facts are obtained immediately from the definition 
of $L_{s}^{\infty}(\X)$ and the triangle inequality, respectively.
\begin{fact}\label{fact1}
Let $F=\chi_{A_{1}}x_{1}+\chi_{A_{2}}x_{2} +\ldots +\chi_{A_{n}}x_{n}$, where $n\in\N$, 
$\{A_{k}|1\leq k\leqq n\}$ is a measurable partition of $[0,1]$ 
and $x_{1},\ldots,x_{n}\in\B_{\X}$. Such simple functions $F$ are dense in $\B_{L_{s}^{\infty}(\X)}$.
\end{fact}

\begin{fact}\label{fact2}
Let $\X$ be a Banach space and $T_{1},...,T_{n}\in \mathcal{G}_{\X},\ n\in \N$. Assume that $x,y,z\in \X$ satisfy
$\dist(y,\conv(\{T_{j}(x)|j=1,...,n\}))=\epsilon>0$ and $||x-z||=\delta> 0$. 
Then 
\[\dist(y, \conv(\{T_{j}(z)|j=1,...,n\}))\leq \epsilon +\delta.\]
\end{fact}

\begin{proof}[Proof of Theorem \ref{th: LPX}]
 
Our strategy for the case $L_{s}^{\infty}(\X)$, which we mainly concentrate on, 
is to show that for a given $x\in \S_{\X}$ the constant function 
$\chi_{[0,1]}x\in L_{s}^{\infty}(\X)$ is a big point and that $\chi_{[0,1]}x\in\overline{\conv}(\mathcal{G}_{L_{s}^{\infty}(\X)}(F))$ 
for each $F\in \S_{L_{s}^{\infty}(\X)}$. This suffices for the claim according to Fact \ref{fact2}.

First we will show that for any $x\in\S_{\X}$ the function $\chi_{[0,1]}x$ is a big point.
By using Fact \ref{fact1} it suffices to show that any simple function $F\in \S_{L_{s}^{\infty}(\X)}$ (as in Fact \ref{fact1}) 
is contained in $\overline{\conv}(\mathcal{G}_{L_{s}^{\infty}(\X)}(\chi_{[0,1]}x))$. Let 
$F=\chi_{A_{1}}y_{1}+\ldots +\chi_{A_{n}}y_{n}$, where $\{A_{k}|1\leq k\leq n\}$ is a measurable partition of $[0,1]$
and $y_{1},\ldots,y_{n}\in \B_{\X}$. Let $\epsilon>0$. 

Since $\X$ is convex-transitive, there are for each $j\in\{1,\ldots,n\}$ finite sets 
$M_{j}\subset\N$, families of isometries $(T_{k}^{(j)})_{k\in M_{j}}\subset\mathcal{G}_{\X}$
and convex weights $(a_{k}^{(j)})_{k\in M_{j}}\in\S_{\ell_{+}^{1}}\cap c_{00}$ such that 
\begin{equation}\label{eq: yj}
\begin{array}{cc}
&\left|\left|\sum_{k\in M_{j}}a_{k}^{(j)}T_{k}^{(j)}(x)-y_{j}\right|\right|<\epsilon\quad \mathrm{for}\ j\in\{1,\ldots,n\}.
\end{array}
\end{equation}
 
Observe that, given $(i_{j})_{j}\in \prod_{j=1}^{n}M_{j}$, one obtains rotations $\boldsymbol{R}\in\mathcal{G}_{L_{s}^{\infty}(\X)}$ by setting
\begin{equation}\label{eq: R*}
\begin{array}{cc}
&\boldsymbol{R}(g)(t)=\sum_{1\leq j\leq n}\chi_{A_{j}}(t)T_{i_{j}}^{(j)}(g(t))\ \mathrm{for\ a.e.}\ t\in [0,1]. 
\end{array}
\end{equation} 
Next we will take a convex combination of this type of rotations. 
By induction on the number of terms in the product one can check that 
\begin{equation}\label{eq: a=1}
\begin{array}{ll}
&\sum_{(i_{1},i_{2},\ldots,i_{n})}\prod_{j=1}^{n}a_{i_{j}}^{(j)}=1,
\end{array}
\end{equation} 
where the summation is taken over all the combinations $(i_{1},i_{2},\ldots,i_{n})\in\prod_{j=1}^{n}M_{j}$. 
One can also verify by induction on the number of terms in the product that
\begin{equation}\label{eq: sumi1}
\begin{array}{rr}
&\sum_{(i_{1},i_{2},\ldots,i_{n})}\left(\prod_{l=1}^{n}a_{i_{l}}^{(l)}\right)\sum_{l=1}^{n}\chi_{A_{l}}(t)T_{i_{l}}^{(l)}=\sum_{l=1}^{n}
\sum_{i_{l}\in M_{l}}a_{i_{l}}^{(l)}\chi_{A_{l}}(t)T_{i_{l}}^{(l)}
\end{array}
\end{equation}
where the leftmost sum is taken over all combinations $(i_{1},\ldots , i_{n})\in \prod_{j=1}^{n}M_{j}$. Indeed, the case $n=1$ is clear. 
To verify the induction step, observe that
\begin{equation*}
\begin{array}{ll}
 &\sum_{(i_{1},\ldots,i_{m+1})}\left(\prod_{l=1}^{m+1}a_{i_{l}}^{(l)}\right)\sum_{l=1}^{m+1}\chi_{A_{l}}(t)T_{i_{l}}^{(l)}\\
=&\sum_{r\in M_{m+1}}a_{r}^{(m+1)}\sum_{(i_{1},\ldots i_{m})}\left(\prod_{l=1}^{m}a_{i_{l}}^{(l)}\right)
\left(\left(\sum_{l=1}^{m}\chi_{A_{l}}(t)T_{i_{l}}^{(l)}\right)+\chi_{A_{m+1}}(t)T_{r}^{(m+1)}\right)\\
=&\sum_{(i_{1},\ldots,i_{m})}\left(\prod_{l=1}^{m}a_{i_{l}}^{(l)}\right)\sum_{l=1}^{m}\chi_{A_{l}}(t)T_{i_{l}}^{(l)}
+\sum_{r\in M_{m+1}}a_{r}^{(m+1)}\chi_{A_{m+1}}(t)T_{r}^{(m+1)}\\
=&\sum_{l=1}^{m}\sum_{i_{l}\in M_{l}}a_{i_{l}}^{(l)}\chi_{A_{l}}(t)T_{i_{l}}^{(l)}+\sum_{r\in M_{m+1}}a_{r}^{m+1}\chi_{A_{m+1}}(t)T_{r}^{(r)}
\end{array}
\end{equation*}
\normalsize
for $m\leq n-1$. Above we used that 
$\sum_{r\in M_{m+1}}a_{r}^{(m+1)}=1=\sum_{(i_{1},\ldots,i_{m})}\prod_{l=1}^{m}a_{i_{l}}^{(l)}$ 
and the second equality is provided by the induction hypothesis.

Hence we obtain by \eqref{eq: sumi1}, \eqref{eq: yj} and the definition of $F$ that
\begin{equation*}
\begin{array}{rl}
   &\left|\left|\sum_{(i_{1},\ldots,i_{n})}\left(\prod_{l=1}^{n}a_{i_{l}}^{(l)}\right)\sum_{l=1}^{n}\chi_{A_{l}}(\cdot)T_{i_{l}}^{(l)}(x)-F(\cdot)\right|\right|_{L_{s}^{\infty}(\X)}\\
  =&\left|\left|\sum_{l=1}^{n}\sum_{i_{l}\in M_{l}}a_{i_{l}}^{(l)}\chi_{A_{l}}(\cdot)T_{i_{l}}^{(l)}(x)-\sum_{l=1}^{n}\sum_{i_{l}\in M_{l}}a_{i_{l}}^{(l)}\chi_{A_{l}}(\cdot)F(\cdot)\right|\right|_{L_{s}^{\infty}(\X)}\leq\epsilon.
\end{array}
\end{equation*}
Note that according to \eqref{eq: R*} and \eqref{eq: a=1} the map 
\begin{equation*}
\begin{array}{cc}
&g\mapsto\sum_{(i_{1},\ldots,i_{n})}\left(\prod_{l=1}^{n}a_{i_{l}}^{(l)}\right)\sum_{l=1}^{n}\chi_{A_{l}}T_{i_{l}}^{(l)}\circ g 
\end{array}
\end{equation*}
above belongs to $\conv(\mathcal{G}_{L_{s}^{\infty}(\X)})$. 
We conclude that $F\in\overline{\conv}(\mathcal{G}_{L_{s}^{\infty}(\X)}(\chi_{[0,1]}x))$, since $\epsilon>0$ was arbitrary. 

Let us check next that for each $F\in \S_{L_{s}^{\infty}(\X)}$ and $x\in\S_{\X}$ it holds that 
$\chi_{[0,1]}x\in \overline{\conv}(\mathcal{G}_{L_{s}^{\infty}(\X)}(F))$. 
To achieve this we employ a kind of sliding hump argument. Fix $F\in \S_{L_{s}^{\infty}(\X)}$.
Since $||F||=1$, one can find by Fact \ref{fact1} a sequence $(x_{n})_{n}\subset\S_{\X}$ and a sequence  
$(B_{n})_{n}$ of measurable (but not necessarily pairwise disjoint) subsets of $[0,1]$ such that 
\begin{equation}\label{eq: supAFx}
\sup_{t\in B_{n}}\ ||F(t)-x_{n}||_{\X}<2^{-(n+1)} 
\end{equation}
and $0<m(B_{n})<\frac{1}{2}$ for all $n$. Put $\Delta_{n}=[1-2^{-n},1-2^{-(n+1)}]$ for $n\in\N$.
We claim that there is a sequence $g_{n}\colon [0,1]\rightarrow [0,1]$ of measurable mappings such that 
\begin{equation}
g_{n}(B_{n})\ae [0,1]\setminus \Delta_{n}\ \mathrm{and}\ g_{n}(B_{n}^{c})\ae \Delta_{n}
\end{equation}
and 
\begin{equation}\label{eq: mequiv}
\mathrm{the\ measure}\ \mu_{n}(\cdot)\stackrel{\cdot}{=}m(g_{n}(\cdot))\colon\Sigma\rightarrow\R\ \mathrm{is\ equivalent\ to}\ m
\end{equation}
for each $n\in \N$. Indeed, let us define two sequences of auxiliary mappings. Define $e_{n}\colon [0,1]\rightarrow [0,1]$ by
\begin{displaymath}
e_{n}(t)=\left\{ \begin{array}{ll}
(1-2^{-(n+1)})\frac{m(B_{n}\cap [0,t])}{m(B_{n})}& \mathrm{if}\ t\in B_{n}\\
(1-2^{-(n+1)})+2^{-(n+1)}\frac{m([0,t]\setminus B_{n})}{m([0,1]\setminus B_{n})}& \mathrm{if}\ t\in [0,1]\setminus B_{n}
\end{array}
\right.
\end{displaymath} 
and $h_{n}\colon [0,1]\rightarrow [0,1]$ by
\begin{displaymath}
h_{n}(t)=\left\{ \begin{array}{ll}
t& \mathrm{if}\ t\in [0,1-2^{-n}]\\
(1-2^{-(n+1)})-(t-(1-2^{-(n+1)}))& \mathrm{if}\ t\in (1-2^{-n},1]
\end{array}
\right.
\end{displaymath}
for $n\in\N$. Then $g_{n}=h_{n}\circ e_{n}$ is the desired mapping for $n\in\N$. Above 
$g_{n}(B_{n})\ae [0,1]\setminus \Delta_{n}$, whose complement $\Delta_{n}$ is the support of the sliding hump.

Next we will apply some observations which appear e.g. in \cite{Greim_Lp} and \cite{Greim_Linfty}. 
Denote by $\Sigma\setminus_{m}$ the quotient sigma algebra of Lebesgue measurable sets on $[0,1]$ formed by identifying the 
$m$-null sets with $\emptyset$. Note that \eqref{eq: mequiv} gives in particular that the map
$\phi_{n}\colon\Sigma\setminus_{m}\rightarrow\Sigma\setminus_{m}$ determined by $\phi_{n}(A)\ae g_{n}(A)$ for $A\in \Sigma$
is a Boolean isomorphism for each $n\in\N$. According to the convex-transitivity of $\X$, there are for each $n\in\N$ a finite set 
$K_{n}\subset\N$, a set of rotations $(T_{i}^{(n)})_{i\in K_{n}}\subset \mathcal{G}_{\X}$ and convex weights 
$(d_{i}^{(n)})_{i\in K_{n}}\in\S_{\ell_{+}^{1}}\cap c_{00}$ 
such that
\begin{equation}\label{eq: xsumd}
\begin{array}{cc}
&\left|\left|x-\sum_{i\in K_{n}}d_{i}^{(n)}T_{i}^{(n)}(x_{n})\right|\right|_{\X}<\frac{1}{2^{n+1}}\quad \mathrm{for\ all}\ n\in\N. 
\end{array}
\end{equation}

Recall that $\psi\colon \Sigma\setminus_{m}\rightarrow \Sigma\setminus_{m}$ is a measure-preserving Boolean isomorphism
if and only if there exists a bijection $\alpha\colon [0,1]\rightarrow [0,1]$ such that $\alpha,\alpha^{-1}$ are measurable and
$\alpha(A)\ae \psi(A)$ for $A\in \Sigma$ (see \cite[p.582-584]{Neu} and \cite[p.340]{HalmosNeumann}). 
Observe that for $n\in\N$ the mapping $\alpha_{n}\colon [0,1]\rightarrow [0,1]$ given by
\begin{displaymath}
\alpha_{n}(t)=\left\{ \begin{array}{ll}
m(B_{n}\cap [0,t])& \mathrm{if}\ t\in B_{n}\\
m(B_{n})+m([0,t]\setminus B_{n})& \mathrm{if}\ t\in [0,1]\setminus B_{n}
\end{array}
\right.
\end{displaymath}
preserves the measure of open intervals of $[0,1]$. Actually, by applying Borel sets we obtain that the mappings $\alpha_{n}$ preserve the 
measure of all measurable sets. For $n\in\N$ let $\psi_{n}$ be the measure-preserving Boolean isomorphism induced by $\alpha_{n}$.
Note that the mappings $g_{n}$ are obtained from $\alpha_{n}$ by simple bijective transformations.
The Boolean isomorphisms $\phi_{n}$ are obtained similarly from $\psi_{n}$.
It follows that for $n\in\N$ there is a measurable bijection $\hat{g}_{n}\colon [0,1]\rightarrow [0,1]$ such that 
$\hat{g}_{n}\circ g_{n}(t)=t$ for $m$-a.e. $t\in [0,1]$ and $\hat{g}_{n}(A)\ae \phi_{n}^{-1}(A)$ for $A\in\Sigma$.
Define mappings $S_{i}^{(n)}\colon L_{s}^{\infty}(\X)\rightarrow L_{s}^{\infty}(\X)$ for $n\in\N$ and $i\in K_{n}$ by
\[S_{i}^{(n)}(F)(t)=T_{i}^{(n)}(F(\hat{g}_{n}(t)))\quad \mathrm{for\ a.e.}\ t\in[0,1],\ F\in L_{s}^{\infty}(\X).\] 
By \eqref{eq: mequiv} we get that $S_{i}^{(n)}\in\mathcal{G}_{L_{s}^{\infty}(\X)}$ (see also \cite[p.467-468]{Greim_Linfty}).

The function $\chi_{[0,1]}x$ can be approximated by convex combinations coming from 
$\conv(\mathcal{G}_{L_{s}^{\infty}(\X)}(F))$ as follows:
\begin{equation}\label{eq: conclusion}
\begin{array}{cc}
&\left|\left|\chi_{[0,1]}x-\frac{1}{k}\sum_{n=1}^{k}\sum_{i\in K_{n}}d_{i}^{(n)}S_{i}^{(n)}(F)\right|\right|_{L_{s}^{\infty}(\X)}\leq \frac{2+\sum_{i=1}^{k}2^{-i} }{k}\stackrel{k\rightarrow \infty}{\longrightarrow}0.
\end{array}
\end{equation}
Indeed, for $n\in \N$ and a.e. $t\in [0,1]\setminus \Delta_{n}$ it holds by \eqref{eq: supAFx} and \eqref{eq: xsumd} that 
\begin{equation*}
\begin{array}{cc}
&\left|\left|x-\sum_{i\in K_{n}}d_{i}^{(n)}S_{i}^{(n)}(F)(t)\right|\right|_{\X}\leq \left|\left|x-\sum_{i\in K_{n}}d_{i}^{(n)}T_{i}^{(n)}(x_{n})\right|\right|_{\X}+\frac{1}{2^{n+1}}
\leq \frac{1}{2^{n}}.
\end{array}
\end{equation*}
On the other hand, $||x-\sum_{i\in K_{n}}d_{n}^{(n)}S_{i}^{(n)}(F)(t)||_{\X}\leq 2$ for a.e. $t\in \Delta_{n}$.
In \eqref{eq: conclusion} we apply the fact that $\Delta_{n}$ are pairwise essentially disjoint. 

Consequently, $\chi_{[0,1]}x\in\overline{\conv}(\mathcal{G}_{L_{s}^{\infty}(\X)}(F))$ for all $F\in \S_{L_{s}^{\infty}(\X)}$ 
by \eqref{eq: conclusion}. We conclude that $L_{s}^{\infty}(\X)$ is convex-transitive.

The case $1\leq p<\infty$ is a straightforward modification of the proof of \cite[Thm. 2.1]{GJK}, where one replaces
$U_{i}x_{i}$ by suitable elements of $\conv(\mathcal{G}_{\X}(x_{i}))$ for each $i$. We will not reproduce the 
details here.
\end{proof}

\subsection{Convex-transitivity of $C_{0}(L,\X)$}
If $L$ is a locally compact space and $\X$ is a Banach space, then $C_{0}(L,\X)$ consists of
the $\X$-valued continuous functions on $L$ tending to $0$ at infinity. This is a Banach space 
endowed with the norm $||f||=\sup_{t\in L}||f(t)||_{\X}$.

Let us recall that for the Cantor set $K$ the space $C(K)$ is convex-transitive but not almost transitive
(see \cite{Rol}). There exist plenty of convex-transitive spaces of the type $C_{0}(L)$ (see e.g. \cite{Wo}).

The following problem about the almost transitivity of $C_{0}(L)$ spaces has attracted attention: 
Wood \cite{Wo} conjectured that if $L$ is a locally compact Hausdorff space and $C_{0}(L)$  is an almost transitive space 
over the scalars $\mathbb{K}$, then $L$ is a singleton. This problem was answered by Greim and Rajagopalan in \cite{GR} 
in the case $\mathbb{K}=\mathbb{R}$ in the \emph{positive} and was recently \emph{refuted} independently by Kawamura in \cite{Kawamura} 
and by Rambla in \cite{Rambla} in the case $\mathbb{K}=\mathbb{C}$. Theorem \ref{thm: singleton} below extends the work in \cite{GR} 
(and here $\mathbb{K}=\R$).

We will require the following generalization of Stone's theorem due to Jerison (see e.g. \cite[p.145]{Beh}):\\
\textit{Let $(L,\tau)$ be a locally compact Hausdorff space and let $\Y$ be a strictly convex Banach space. 
Then a map $T\colon C_{0}(L,\Y)\rightarrow C_{0}(L,\Y)$ is a rotation if and only if $T$ has the form
\begin{equation}\label{eq: jerison}
T(F)(t)=\sigma(t)(F(h(t)))\quad \mathrm{for}\ F\in C_{0}(L,\Y);\ t\in L,
\end{equation}
where $\sigma\colon L\rightarrow \mathcal{G}_{\Y}$ is $(\tau,\mathrm{SOT})$-continuous and $h\colon L\rightarrow L$ is a homeomorphism.}
\begin{theorem}\label{thm: singleton}
Let $(L,\tau)$ be locally compact Hausdorff space that is locally connected at some point $x\in L$. 
Suppose that $\X$ is a non-trivial strictly convex Banach space. 
\begin{enumerate}
\item[(i)]{If $C_{0}(L,\X)$ is almost transitive, then $L$ is a singleton.}
\item[(ii)]{If $C_{0}(L,\X)$ is convex-transitive, then $L$ is connected.}
\end{enumerate}
\begin{proof}
Let us first make some preparations towards both these claims. Let $x,y\in L$ be such that $x\neq y$ and $L$ is locally connected at $x$. 
Since $L$ is completely regular and Hausdorff, there are open neighbourhoods $V_{x}$ and $V_{y}$ of $x$ and $y$, respectively,
such that $\overline{V}_{x}^{\tau}\cap \overline{V}_{y}^{\tau}=\emptyset$. According to the complete regularity of 
$L$, there is a continuous map $f_{0}\colon L\rightarrow [0,1]$ such that $f_{0}(L\setminus V_{x})=\{0\}$ and $f_{0}(x)=1$.  
Since $L$ is locally connected at $x$ and $f_{0}^{-1}((\frac{1}{2},1])$ is an open neighbourhood
of $x$, there is an open connected neighbourhood $U_{x}\subset f_{0}^{-1}((\frac{1}{2},1])$ of $x$. 
Again, by the complete regularity there is a positive continuous map $e\colon L\rightarrow [0,1]$ such that
$e(L\setminus U_{x})=\{0\}$ and $e(x)=1$. Define $f=\min(f_{0},\frac{1}{2})+\max(e\cdot(f_{0}-\frac{1}{2}),0)$. Note that 
$f\in C_{0}(L)$, its range is in $[0,1]$, $f(L\setminus U_{x})\subset [0,\frac{1}{2}]$ and $f(x)=1$. 
Let $g\colon L\rightarrow [0,1]$ be a positive map in $C_{0}(L)$ such that $g(L\setminus V_{y})=\{0\}$ and $g(y)=1$.

Fix $w\in \S_{\X}$. Note that $f\otimes w, (f+g)\otimes w \in \S_{C_{0}(L,\X)}$ by the construction of $f$ and $g$. 

\textit{Claim} (i). Suppose that $x,y\in L$ are as above. Then $C_{0}(L,\X)$ is not almost transitive. 

Indeed, assume to the contrary that $C_{0}(L,\X)$ is almost transitive. Hence there must be a rotation 
$T\colon C_{0}(L,\X)\rightarrow C_{0}(L,\X)$ such that 
\begin{equation}\label{eq: epsilon10}
|\ ||T(f\otimes w)||-||(f+g)\otimes w||\ |\leq ||T(f\otimes w)-(f+g)\otimes w||<\frac{1}{10}.
\end{equation} 
By \eqref{eq: jerison} we know that $T(Q\otimes w)(t)=\sigma(t)(Q(h(t))w)$ for any $Q\in C_{0}(L)$, 
where $h\colon L\rightarrow L$ is a homeomorphism and $\sigma\colon L\rightarrow \mathcal{G}_{\X}$. 
In particular 
\begin{equation}\label{eq: | |}
||T(Q\otimes w)(t)||_{\X}=||Q(h(t))w||_{\X}=|Q(h(t))|\ \mathrm{for\ each}\ t\in L,
\end{equation}
since $||w||_{\X}=1$.
Thus by \eqref{eq: epsilon10}, \eqref{eq: | |} and the positivity of $f$ and $g$ we obtain that
\begin{equation}\label{eq: lip10}
f(h(t))-\frac{1}{10}< (f+g)(t)< f(h(t))+\frac{1}{10}\quad \mathrm{for}\ t\in L.
\end{equation}
Note that since $f(x)=g(y)=1$, and $f(y)=g(x)=0$, it holds that $f(h(x)),f(h(y))\geq \frac{9}{10}$. 
Since $f(L\setminus U_{x})\subset [0,\frac{1}{2}]$, we obtain that $h(x),h(y)\in U_{x}$ and hence $x,y\in h^{-1}(U_{x})$.
By \eqref{eq: lip10} and the selection of $U_{x}$ we have that 
\[(f+g)(h^{-1}(t))>f(t)-\frac{1}{10}> \frac{1}{2}-\frac{1}{10}=\frac{2}{5}\]
for all $t\in U_{x}$. Since $f$ and $g$ are disjointly supported, we get that
\[x,y\in h^{-1}(U_{x})\subset f^{-1}((2/5,1])\ \cup\ g^{-1}((2/5,1]).\]
Moreover, there exists a separation between the open sets $f^{-1}((2/5,1])\ni x$ and 
$g^{-1}((2/5,1])\ni y$. Thus $h^{-1}(U_{x})$ is not connected and hence $U_{x}$ is not connected; a contradiction.

\textit{Claim} (ii). Suppose that $x,y\in L$ are as in the beginning of the proof. Since $y$ was arbitrary, it suffices to
show that $x$ and $y$ belong to the same connected component. 

Let $w^{\ast}\in S_{\X^{\ast}}$ be a support functional for $w\in \S_{\X}$.
Define $\boldsymbol{F}\in \S_{C_{0}(L,\X)^{\ast}}$ by $\boldsymbol{F}(u)=\frac{w^{\ast}(u(x)-u(y))}{2}$ for each $u\in C_{0}(L,\X)$.
Note that $\boldsymbol{F}$ is a support functional for $(f-g)\otimes w\in \S_{C_{0}(L,\X)}$ as 
\[\boldsymbol{F}((f-g)\otimes w)=\frac{f(x)-g(x)-(f(y)-g(y))}{2}w^{\ast}(w)=\frac{f(x)+g(y)}{2}=1.\] 
Since $C_{0}(L,\X)$ was assumed to be convex-transitive we obtain that 
\[\sup \boldsymbol{F}(\overline{\conv}(\mathcal{G}_{C_{0}(L,\X)}(f\otimes w)))=\boldsymbol{F}((f-g)\otimes w)=1.\]
In particular we may select by the linearity of $\boldsymbol{F}$ an element
$G\in \mathcal{G}_{C_{0}(L,\X)}(f\otimes w)$ such that $\boldsymbol{F}(G)>\frac{9}{10}$. This means that 
$w^{\ast}(G(x)),-w^{\ast}(G(y))\in (\frac{8}{10},1]$. Hence 
\begin{equation}\label{eq: gxz}
||G(x)||_{\X},||G(y)||_{\X}\in (4/5,1], 
\end{equation}
since $||w^{\ast}||_{\X^{\ast}}=1$. 

Again, by \eqref{eq: jerison} we may write $G(t)=\sigma(t)(f(h(t))w)$. Note that 
$||G(t)||_{\X}=|f(h(t))|$ for all $t\in L$. Thus $||G(t)||_{\X}\leq \frac{1}{2}$ for all 
$t\in L\setminus h^{-1}(U_{x})$ and hence $x,y \in h^{-1}(U_{x})$ by \eqref{eq: gxz}. 
Clearly $x$ and $y$ belong to the same connected component as $U_{x}$ is connected.
\end{proof}
\end{theorem}

\begin{theorem}
Let $(L,\tau)$ be a locally compact Hausdorff space. 
\begin{enumerate}
\item[(i)]{If $\X$ is a strictly convex Banach space and $C_{0}(L,\X)$ is convex-transitive, then $\X$ is convex-transitive.} 
\item[(ii)]{If $\H$ is a Hilbert space, $\dim(\H)=\infty$ and $C_{0}(L)$ is convex-transitive, then $C_{0}(L,\H)$ is convex-transitive.}
\end{enumerate}
\end{theorem}
Our argument below relies strongly on the assumption that $\dim(\H)=\infty$. We do not know if this assumption can be removed.
\begin{proof}
(i). Fix $x,y\in\S_{\X}$. Let $F,G\in C_{0}(L,\X)$ be given by $F=e(\cdot)x,\ G=e(\cdot)y$, where $e\in \S_{C_{0}(L)}$.
By the local compactness of $L$ and the continuity of $e$ there is $t_{0}\in L$ such that $|e(t_{0})|=1$. 
Since $C_{0}(L,\X)$ is convex-transitive, there exists a sequence 
$(T_{n})_{n\in\N}\subset \mathcal{G}_{C_{0}(L,\X)}$ such that $G\in \overline{\conv}((T_{n}(F))_{n\in\N})$. Since $\X$ is strictly convex, 
representation \eqref{eq: jerison} gives that $T_{n}(F)(t_{0})=\sigma_{n}(t_{0})(e(h_{n}(t_{0}))x)$ for suitable
sequences $(\sigma_{n})$ and $(h_{n})$. Recall that here $\sigma_{n}\colon L\rightarrow \mathcal{G}_{\X}$ are continuous and 
$h_{n}\colon L\rightarrow L$ are homeomorphisms for $n\in\N$.
Since $G(t_{0})\in \{\pm y\}$ we obtain for $a_{n}\stackrel{\cdot}{=}e(h_{n}(t_{0}))\in [-1,1]$ 
that 
\begin{equation}
\begin{array}{rl}
y\in &\overline{\conv}(\{\sigma_{n}(t_{0})(a_{n}x)|n\in\N\})\subset\overline{\conv}(\conv(\{\sigma_{n}(t_{0})(\pm x)|n\in\N\}))\\
    =&\overline{\conv}(\{\pm\sigma_{n}(t_{0})(x)|n\in\N\})\subset\overline{\conv}(\mathcal{G}_{\X}(x)).
\end{array}
\end{equation}
Hence $\X$ is convex-transitive.

(ii). Let $\H$ be an infinite-dimensional Hilbert space and let $\{e_{\gamma}\}_{\gamma\in\Gamma}$ be an orthonormal basis of $\H$.
Fix $\gamma_{0}\in \Gamma$, $0<\epsilon<1$ and $F_{0},G_{0}\in \S_{C_{0}(L,\H)}$.
Define $\mathrm{r}\colon \H\rightarrow \S_{\H}\cup\{0\}$ by $\mathrm{r}(z)=\frac{z}{||z||}$ for $z\neq 0$ and $\mathrm{r}(0)=0$.
Recall that $\mathrm{r}$ is continuous away from $0$.

We will make crucial use of the fact that there exists a $(||\cdot||,\mathrm{SOT})$ -continuous map
$A\colon \S_{\H}\setminus \{e_{\gamma_{0}}\}\rightarrow \mathcal{G}_{\H}$ such that 
\begin{equation}\label{eq: Ax}
A(x)(e_{\gamma_{0}})=x\ \mathrm{for\ all}\ x\in\S_{\H}\setminus\{e_{\gamma_{0}}\}, 
\end{equation}
see \cite[Prop 5.4]{AR}. We will advance in two steps to complete the proof. 
The first, rather technical step is needed in order to apply the above fact \eqref{eq: Ax}. This step uses the assumption 
that $\dim(\H)=\infty$. The second step involves applying the fact \eqref{eq: jerison}.

\textit{Step 1.} Note that the subset $\{t\in L:\ ||F_{0}(t)||\geq \frac{\epsilon}{8}\}\subset L$ is compact, and its image under the 
continuous mapping $F_{0}$ is also compact. Since $\dim(\H)=\infty$ one can find by approximation an index 
$\mu\in \Gamma,\ \mu\neq \gamma_{0},$ such that $||P_{\mu}\circ F_{0}(t)||<\frac{\epsilon}{8}$ 
for $t\in L$ such that $||F_{0}(t)||\geq \frac{\epsilon}{8}$. Above $P_{\mu}\colon \H\rightarrow [e_{\mu}]$ 
is the orthogonal projection. It follows that $||P_{\mu}\circ F_{0}(t)||<\frac{\epsilon}{8}$ for all $t\in L$.

Towards an application of \eqref{eq: Ax} we define an auxiliary mapping $F\in C_{0}(L,\H)$ by perturbing $F_{0}$ as follows:
\begin{displaymath}
F(t)=\left\{ \begin{array}{ll}
||F_{0}(t)||\mathrm{r}\big(F_{0}(t)+\frac{\epsilon}{8}e_{\mu}\big)& \mathrm{if}\ 
||F_{0}(t)||\geq\frac{\epsilon}{2}\\
||F_{0}(t)||\mathrm{r}\big(||F_{0}(t)||\mathrm{r}\big(F_{0}(t)+\frac{\epsilon}{8}e_{\mu}\big)+(\frac{\epsilon}{8}-\frac{1}{4}||F_{0}(t)||)e_{\mu}\big)& \mathrm{if}\ ||F_{0}(t)||<\frac{\epsilon}{2}
\end{array}
\right.
\end{displaymath} 
Put $\mathrm{r}_{F}(t)=\mathrm{r}(F(t))$ for $t\in L$ such that $F(t)\neq 0$ and $\mathrm{r}_{F}(t)=e_{\mu}$ otherwise.
Then $F$ and $\mathrm{r}_{F}$ satisfy the following conditions:
\begin{enumerate}
\item[(a)]{$||F(t)||=||F_{0}(t)||$ for $t\in L$.}
\item[(b)]{$||F_{0}-F||<\epsilon$.}
\item[(c)]{$F(L)\cap [e_{\gamma_{0}}]\subset \{0\}$.}
\item[(d)]{$\underset{||F_{0}(t)||\rightarrow 0}{\lim}\mathrm{r}_{F}(t)=e_{\mu}$.}
\end{enumerate}
Indeed, conditions (a) and (d) are immediate. Condition (b) follows by analyzing the upper and lower cases separately. 
In the upper case one applies geometric estimates in the euclidean plane spanned by $F_{0}(t)$ and $e_{\mu}$ in $\H$.
If one writes $||F_{0}(t)||\mathrm{r}(F_{0}(t)+\frac{\epsilon}{8}e_{\mu})=b F_{0}(t)+c e$ for suitable $b,c\in\R$ and $e\bot F_{0}(t)$, 
where $||e||=1$ and $||F_{0}(t)||\geq \frac{\epsilon}{2}$, then $|(||F_{0}(t)||-b)|\leq \frac{\epsilon}{8}$ and 
$|c|\leq \frac{\epsilon}{4}$. The lower case follows immediately by applying the triangle inequality.
Towards condition (c), given $t\in L$, write $P_{\mu}(F_{0}(t)+\frac{\epsilon}{2}e_{\mu})=a e_{\mu}$ for 
suitable $a\in\R$. The fact that $||P_{\mu}\circ F_{0}(t)||<\frac{\epsilon}{8}$ for $t\in L$ yields 
that $a>0$ above.  
Thus $P_{\mu}\mathrm{r}\big(F_{0}(t)+\frac{\epsilon}{8}e_{\mu}\big)\neq 0$ for all $t\in L$ and 
$P_{\mu} \mathrm{r}\big(||F_{0}(t)||\mathrm{r}\big(F_{0}(t)+\frac{\epsilon}{8}e_{\mu}\big)+(\frac{\epsilon}{8}-\frac{1}{4}||F_{0}(t)||)e_{\mu}\big)\neq 0$ 
for $t$ such that $||F_{0}(t)||<\frac{\epsilon}{2}$. Hence (c) holds, since $e_{\mu}\bot e_{\gamma_{0}}$.

Similarly, by using a suitable $\nu\in\Gamma\setminus\{e_{\gamma_{0}}\}$, one can define $G\in C_{0}(L,\H)$ and $\mathrm{r}_{G}$, 
which correspond to $G_{0}$ and have analogous properties. 
Observe that $\mathrm{r}_{F}$ and $\mathrm{r}_{G}$ are continuous mappings $L\rightarrow \S_{\H}$ by condition (d) 
and the continuity of $F$ and $G$. 

\textit{Step 2.} We claim that the map 
\[E(t)\mapsto A(\mathrm{r}_{G}(t))\circ A(\mathrm{r}_{F}(t)))^{-1}(E(t)),\quad t\in L,\ E\in C_{0}(L,\H)\]
defines an element in $\mathcal{G}_{C_{0}(L,\H)}$. Indeed, recall that $A$ is $(||\cdot||,\mathrm{SOT})$-continuous and 
$(\mathcal{G}_{\H},\mathrm{SOT})$ is a topological group, which means that the group operation and the taking of the inverse 
are continuous. Thus the composition 
$L\rightarrow \S_{\H}\times \S_{\H}\rightarrow \mathcal{G}_{\H}$ given by 
\[t\mapsto (\mathrm{r}_{G}(t),\mathrm{r}_{F}(t))\mapsto A(\mathrm{r}_{G}(t))\circ (A(\mathrm{r}_{F}(t)))^{-1}\] 
is $(\tau,\mathrm{SOT})$-continuous. Clearly 
$A(\mathrm{r}_{G}(t))\circ (A(\mathrm{r}_{F}(t)))^{-1}\in \mathcal{G}_{\H}$ for $t\in L$. By \eqref{eq: Ax} we obtain that
\begin{equation}\label{eq: signGt}
\mathrm{r}_{G}(t)=A(\mathrm{r}_{G}(t))(e_{\gamma_{0}})=A(\mathrm{r}_{G}(t))\circ (A(\mathrm{r}_{F}(t)))^{-1}(\mathrm{r}_{F}(t))\quad \mathrm{for}\ t\in L.  
\end{equation}

Consider $||G(\cdot)||,||F(\cdot)||\in C_{0}(L)$. One can find by the convex-transitivity of $C_{0}(L)$ and \eqref{eq: jerison} 
homeomorphisms $h_{n}\colon L\rightarrow L$ and continuous functions $\theta_{n}\colon L\rightarrow \{-1,1\}$ for $n\in\N$ 
such that $||G(\cdot)||\in \overline{\conv}(\{\theta_{n}(\cdot)||F(h_{n}(\cdot))||:\ n\in\N\})$.  

Put $T_{n}(E)(t)=A(\mathrm{r}_{G}(t))\circ (A(\mathrm{r}_{F}(h_{n}(t))))^{-1}(E(t))$ for $n\in\N,\ E\in C_{0}(L,\H)$ and $t\in L$. 
By \eqref{eq: signGt} we obtain that $T_{n}(F\circ h_{n})(t)=||F(h_{n}(t))||\mathrm{r}_{G}(t)$ for all $t\in L$.

Since $||G(\cdot)||\in \overline{\conv}(\{\theta_{n}(\cdot)||F(h_{n}(\cdot))||:\ n\in\N\})\subset C_{0}(L)$, we obtain that 
\begin{eqnarray*}
  G(\cdot)&=&||G(\cdot)||\mathrm{r}_{G}(\cdot)\in \overline{\conv}(\{\theta_{n}(\cdot)||F(h_{n}(\cdot))||\mathrm{r}_{G}(\cdot):\ n\in\N\})\\ 
          &=&\overline{\conv}(\{(\theta_{n}(\cdot)\I)\circ T_{n}(F(h_{n}(\cdot))):\ n\in\N\})\subset C_{0}(L,\H).
\end{eqnarray*}
Note that according to \eqref{eq: jerison} the map 
$E(\cdot)\mapsto (\theta_{n}(\cdot)\I)\circ T_{n}(E(h_{n}(\cdot)))$ 
defines an element in $\mathcal{G}_{C_{0}(L,\H)}$ for $n\in\N$. Hence $G\in\overline{\conv}(\mathcal{G}_{C_{0}(L,\H)}(F))$. 
Since $F_{0},G_{0}$ and $\epsilon$ were arbitrary, Fact \ref{fact2} and condition (b) yield that $C_{0}(L,\H)$ is convex-transitive.
\end{proof}

\section{Some new examples}
Even though there exist plenty of transitive spaces, very few \emph{concrete} examples of almost transitive and convex-transitive 
(scalar-valued) function spaces are known (see \cite{BR2}). Next we will demonstrate transitivity conditions for a few fairly concrete 
spaces by direct proofs. These new examples are obtained by performing natural operations on classical function spaces.

In this section we denote 
\[\ell^{\infty}_{\lambda}(\kappa)=\{x\in \ell^{\infty}(\kappa):\ |\mathrm{supp}(x)|\leq \lambda\}\subset \ell^{\infty}(\kappa),\]
where $\kappa,\ \lambda$ are infinite cardinals and $\lambda\leq \kappa$. This clearly defines a closed subspace of
$\ell^{\infty}(\kappa)$. Recall that $\kappa^{+}$ stands for the successor cardinal of $\kappa$, i.e. the least cardinal strictly
larger than $\kappa$, and $\kappa^{++}=(\kappa^{+})^{+}$. 

In \cite{asytrans} the special role of $1$-codimensional subspaces of $L^{p}$ in connection with the rotations is discussed.
There it is also verified that the subspace 
\[M^{1}\stackrel{\cdot}{=}\left\{x\in L^{1}\Big{|}\int_{0}^{1}x(t)\ \mathrm{d}t=0\right\}\subset L^{1}\] 
appearing in the following result is neither $1$-complemented, nor isometric to $L^{1}$.
\begin{theorem}
The subspace $M^{1}\subset L^{1}$ is almost transitive.
\end{theorem}
\begin{proof}
For each $f\in L^{1}$ we denote $f_{+}(t)=\max(f(t),0)$ and $f_{-}(t)=\min(f(t),0)$, 
so that $f=f_{+}+f_{-}$. Fix $x,y\in \S_{M^{1}}$ and $\epsilon>0$. 
Let $A\cup B=[0,1]$ be a measurable partition such that $A=\mathrm{supp}(x_{+})$.
Similarly, let $C\cup D=[0,1]$ be a measurable partition such that $C=\mathrm{supp}(y_{+})$.
Observe that $||x_{+}||=||x_{-}||=||y_{+}||=||y_{-}||=\frac{1}{2}$. 

Note that if $E,F\subset [0,1]$ are measurable sets with strictly positive measure and $R\colon L^{1}(E)\rightarrow L^{1}(F)$ 
is a positive linear isometry onto, then 
$\int_{E}f(t)\ \mathrm{d}t=\int_{F}R(f)(t)\ \mathrm{d}t$ for $f\in L^{1}$. Indeed, 
\begin{equation*}
\begin{array}{rl}
&\int_{E}f(t)\ \mathrm{d}t=\int_{E}f_{+}(t)\ \mathrm{d}t + \int_{E}f_{-}(t)\ \mathrm{d}t=||f_{+}||_{L^{1}(E)}-||f_{-}||_{L^{1}(E)}\\
=&||R(f_{+})||_{L^{1}(F)}-||R(f_{-})||_{L^{1}(F)}=\int_{F}R(f_{+})(t)\ \mathrm{d}t-\int_{F}-R(f_{-})(t)\ \mathrm{d}t\\
=&\int_{F}R(f)(t)\ \mathrm{d}t,
\end{array}
\end{equation*}
since $f=f_{+}+f_{-}$, $R(f_{+})\geq 0$ and $R(f_{-})\leq 0$.

Let $R_{1}\colon L^{1}(A)\rightarrow L^{1}(C)$ and $R_{2}\colon L^{1}(B)\rightarrow L^{1}(D)$ be positive linear isometries onto.
Then $\int_{A}x_{+}(t)\ \mathrm{d}t=\int_{C}R_{1}(x_{+}|_{A})(t)\ \mathrm{d}t$ and 
$\int_{B}x_{-}(t)\ \mathrm{d}t=\int_{D}R_{2}(x_{-}|_{B})(t)\ \mathrm{d}t$ by the above observation. 
Since both $R_{1}(x_{+}|_{A}),y_{+}|_{C}\in L^{1}(C)$ are 
positive and have norm $\frac{1}{2}$, it follows by the standard argument, which is used to prove that $L^{1}$ is almost transitive
(see e.g. \cite[Thm. 3.1]{Lamperti}), that there is a positive rotation $T_{1}\colon L^{1}(C)\rightarrow L^{1}(C)$ such that 
$||T_{1}R_{1}(x_{+}|_{A})-y_{+}|_{C}||<\frac{\epsilon}{2}$. 
Similarly there is a positive rotation $T_{2}\colon L^{1}(D)\rightarrow L^{1}(D)$ such that
$||T_{2}R_{2}(x_{-}|_{B})-y_{-}|_{D}||<\frac{\epsilon}{2}$. Now 
\[T_{1}R_{1}\oplus T_{2}R_{2}\colon L^{1}=L^{1}(A)\oplus_{1}L^{1}(B)\rightarrow L^{1}(C)\oplus_{1}L^{1}(D)=L^{1}\] 
is a positive linear isometry onto. Hence $\int_{[0,1]}f(t)\ \mathrm{d}t=\int_{[0,1]}R(f)(t)\ \mathrm{d}t$ for each
$f\in L^{1}$, so that $(T_{1}R_{1}\oplus T_{2}R_{2})|_{M^{1}}\colon M^{1}\rightarrow M^{1}$ is a linear isometry onto.
Note that $||(T_{1}R_{1}\oplus T_{2}R_{2})x-y||<\epsilon$. Since $x,y$ and $\epsilon$ were arbitrary, we conclude that
$M^{1}$ is almost transitive.
\end{proof}

We refer to \cite{Heinrich} for background in ultraproducts in Banach space theory. 
Recall that an ultrafilter $\mathcal{U}$ is non-principal if it does not contain any singletons.
The writer thanks \AA. Hirvonen and T. Hyttinen for providing help in proving the following fact. 
\begin{lemma}(Homogeneity)\label{UFlemma} 
Let $\mathcal{U}$ be a non-principal ultrafilter on $\N$ and $A,B\subset \N$ be infinite sets such that $A,B\notin \mathcal{U}$.
Then there is a permutation $h\colon \N \rightarrow \N$ such that
\begin{enumerate}
\item[(i)]{$h(A)=B$.}
\item[(ii)]{$h(U)\in \mathcal{U}$ if and only if $U\in\mathcal{U}$.}
\end{enumerate}
\end{lemma} 
\begin{proof}
Observe that $\N\setminus (A\cup B)\in \mathcal{U}$. Let $L,M\subset \N\setminus (A\cup B)$ be infinite disjoint subsets
such that $L\cup M=\N\setminus (A\cup B)$. Without loss of generality assume that $M\in \mathcal{U}$ and 
$L\notin \mathcal{U}$.

Since $A$ and $L$ are infinite and disjoint, there is a bijection $g\colon A\cup L\rightarrow A\cup L$ such that 
$g(A)=L$ and $g(L)=A$. Define $\tilde{g}\colon \N\rightarrow \N$ by $\tilde{g}|_{A\cup L}=g$ and 
$\tilde{g}|_{\N\setminus (A\cup L)}=\mathrm{id}|_{\N\setminus (A\cup L)}$. Note that $V\in\mathcal{U}$
if and only if $V\setminus (A\cup L)\in \mathcal{U}$. Hence $\tilde{g}(U)\in \mathcal{U}$ if and only if 
$U\in\mathcal{U}$. 

Similarly, one can construct a bijection $\tilde{f}\colon \N\rightarrow \N$ such that $\tilde{f}(L)=B$, 
$\tilde{f}(B)=L$ and $\tilde{f}(U)\in \mathcal{U}$ if and only if $U\in\mathcal{U}$. Observe that 
$h=\tilde{f}\circ\tilde{g}\colon \N\rightarrow \N$ is the claimed permutation.
\end{proof}

\begin{fact}\label{Blinfty}
Let $\Gamma$ be a non-empty set. Then $\B_{\ell^{\infty}(\Gamma)}=\overline{\conv}(\{-1,1\}^{\Gamma})\subset\ell^{\infty}(\Gamma)$.
\end{fact}
The proof is left as an exercise (see \cite[p.99]{HHZ}). Note that each $x\in \B_{\ell^{\infty}(\Gamma)}$ can be approximated 
by elements of $\{-1,-\frac{n-1}{n},\ldots,\frac{n-1}{n},1\}^{\Gamma}$, where $n\in\N$.
 
\begin{theorem}\label{th prev}
Let $\mathcal{U}$ be a non-principal ultrafilter on $\N$ and  
$$\Y=\{ (x_{n})\in \ell^{\infty}|\ \lim_{\mathcal{U}}x_{n}=0\}.$$ 
Then $\Y\mod c_{0}$ is convex-transitive.
\end{theorem}

Let us make a few comments before the proof. Recall that $\ell^{\infty}$ itself is far away from being convex-transitive
in view of the characterization of the rotations of $\ell^{\infty}$ (see \cite[2.f.14]{LTI}). 
Note that \emph{exactly} one of the sets $\{2n|\ n\in\N\}$ and $\{2n+1|\ n\in\N\}$
is a member of $\mathcal{U}$, say $\{2n+1|n\in\N\}\in\mathcal{U}$. Thus $(x_{n})\mapsto (0,x_{1},0,x_{2},0,x_{3},\ldots)$
defines an isometric embedding $\ell^{\infty}\rightarrow \Y$. In particular, $\Y\mod c_{0}$ is a non-separable space
containing an isometric copy of $\ell^{\infty}\mod c_{0}$.
\begin{proof}
First note that if $x=(x_{n})\in\ell^{\infty}$ is such that $\lim_{\mathcal{U}}|x_{n}|=c>0$, then 
$\dist(x,\Y)\geq c>0$ because $||x-y||\geq c$ for each $y\in\Y$. Thus $\Y\subset \ell^{\infty}$ is a closed subspace.

Recall (see \cite[2.f.14]{LTI}) that $T\in \mathcal{G}_{\ell^{\infty}}$ if and only if
\begin{equation}\label{eq: Tsum} 
T((x_{n})_{n})=(\theta(n)x_{\pi(n)})_{n},
\end{equation}
where $\pi\colon \N\rightarrow\N$ is a bijection and $\theta\colon \N\rightarrow \{-1,1\}$. 
Also note that such an isometry $T$ restricted to $c_{0}$ is a member of $\mathcal{G}_{c_{0}}$.

Next we will check that if $T\in\mathcal{G}_{\ell^{\infty}}$, then $\widehat{T}\colon x+c_{0}\mapsto T(x)+c_{0}$, for $x\in\ell^{\infty}$,
defines a rotation $\ell^{\infty}\mod c_{0}\rightarrow \ell^{\infty}\mod c_{0}$. 
Indeed, it is well-known that $\widehat{T}\colon \ell^{\infty}\mod c_{0}\rightarrow \ell^{\infty}\mod c_{0}$ is a linear bijection. 
Moreover, $\inf_{z\in c_{0}}||x-z||=\inf_{z\in c_{0}}||T(x)-T(z)||=\inf_{z\in c_{0}}||T(x)-z||$, so that 
$\widehat{T}\colon \ell^{\infty}\mod c_{0}\rightarrow \ell^{\infty}\mod c_{0}$ is an isometry.

Note that if the permutation $\pi$ above is such that $\pi^{-1}(U)\in \mathcal{U}$ exactly when $U\in\mathcal{U}$, then 
$T|_{\Y}\colon \Y\rightarrow \Y$ is a rotation. To summarize, $\widehat{T}|_{\Y\mod c_{0}}\in\mathcal{G}_{\Y\mod c_{0}}$ for such $\pi$.

Fix $u,v\in \S_{\Y\mod c_{0}}$. If $x,y\in \Y$ are such that $u=x+c_{0}$ and $v=y+c_{0}$, then 
\begin{equation}\label{eq: distsup}
\dist(x,c_{0})=\limsup_{n\rightarrow\infty}|x_{n}|=1=\dist(y,c_{0})=\limsup_{n\rightarrow\infty}|y_{n}|, 
\end{equation}
since $u,v\in \S_{\Y\mod c_{0}}$. Hence we may pick $x,y\in \S_{\Y}$ such that $u=x+c_{0}$ and $v=y+c_{0}$.    

Fix $k\in\N$. Let $A=\{n\in \N:\ |x_{n}|\geq 1-\frac{1}{k}\}$ and $B=\{n\in\N:\ |y_{n}|\geq \frac{1}{k}\}$. 
Observe that these are \emph{infinite} sets by \eqref{eq: distsup} and that $\N\setminus A,\ \N\setminus B\in \mathcal{U}$. 
Thus, by Lemma \ref{UFlemma} 
there is a bijection $\pi\colon \N\rightarrow \N$ such that $\pi(B)=A$ and $\pi(U)\in\mathcal{U}$ 
if and only if $U\in\mathcal{U}$. Observe that $T\colon \Y\rightarrow \Y;\ (x_{n})_{n\in\N}\mapsto (x_{\pi(n)})_{n\in\N}$
is a rotation. Put $w=(w_{n})_{n\in\N}=(x_{\pi(n)})_{n\in\N}\in\mathcal{G}_{\Y}(x)$. 
By the definition of the sets $A$ and $B$ we obtain that $|w_{n}|\geq 1-\frac{1}{k}$ if and only if $n\in B$.

Fix an auxiliary point $z\in\Y$, which satisfies 
\[||y-z||_{\Y}\leq\frac{1}{k},\ |z_{n}|\leq 1-\frac{1}{k}\ \mathrm{for}\ n\in B\ \mathrm{and}\ z_{n}=0\ \mathrm{for}\ n\in \N\setminus B.\]

Write $\ell^{\infty}=\ell^{\infty}(B)\oplus_{\infty}\ell^{\infty}(\N\setminus B)$. It follows from Fact \ref{Blinfty} that
\[z|_{B}\in\overline{\conv}(\{(\theta(n)|\ \theta(n)=\pm 1\ \mathrm{for}\ n\in B \}).\] 
Since 
\begin{equation}\label{eq: wsignw}
||(w_{n})_{n\in B}-(\sign(w_{n}))_{n\in B}||_{\ell^{\infty}(B)}\leq \frac{1}{k} 
\end{equation}
we obtain by Fact \ref{fact2} that
\begin{equation}\label{eq: distzbconv}
 \dist(z|_{B},\overline{\conv}(\{(\theta(n)w_{n})_{n\in B}|\ \theta(n)=\pm 1\ \mathrm{for}\ n\in B \})\leq \frac{1}{k}.
\end{equation}
It follows that 
\begin{equation}\label{eq: distz}
\dist(z,\overline{\conv}(\{(\theta(n)w_{n})_{n\in\N}|\ \theta(n)=\pm 1\ \mathrm{for}\ n\in\N\}))\leq \frac{1}{k}.
\end{equation}
Indeed, one can use convex combinations 
\begin{equation}\label{eq: 4.3*}
\frac{1}{2}\left(\sum_{i}a_{i}(f_{i}+g)+\sum_{i}a_{i}(f_{i}-g)\right)\in \ell^{\infty},
\end{equation}
where $f_{i}\in \{(\theta(n)w_{n})_{n\in B}|\ \theta(n)=\pm 1\ \mathrm{for}\ n\in B \}$ are obtained from \eqref{eq: distzbconv}
such that $||\sum_{i}a_{i}f_{i}-z|_{B}||_{\ell^{\infty}(B)}\leq \frac{1}{k}$ 
and $g=(w_{n})_{n\in \N\setminus B}\in \ell^{\infty}(\N\setminus B)$. Above we used that 
$z(n)=0$ for $n\in \N\setminus B$. 

Since $||y-z||\leq \frac{1}{k}$ by the selection of $z$ and $(w_{n})=(x_{\pi(n)})\in \mathcal{G}_{\Y}(x)$, Fact \ref{fact2} yields that 
\[\dist(y,\overline{\conv}(\mathcal{G}_{\Y}(x)))\leq \dist(y,\overline{\conv}(\{(\theta(n)x_{\pi(n)})_{n\in\N}|\ \theta\colon \N\rightarrow \{\pm 1\}\ \}))\leq \frac{2}{k}.\]
Since $k$ was arbitrary we obtain that $y\in\overline{\conv}(\mathcal{G}_{\Y}(x))$. Thus 
$v=y+c_{0}\in \overline{\conv}(\mathcal{G}_{\Y\mod c_{0}}(u))$, so that $\Y\mod c_{0}$ is convex-transitive.
\end{proof}
We note that the technique applied above does not give that $\Y$ or $\ell^{\infty}\mod c_{0}$ should be convex-transitive. However,
by applying similar ideas we obtain another fairly concrete example of a convex-transitive space:
\begin{theorem}
The space $\ell^{\infty}_{\kappa^{+}}(\kappa^{++})\mod \ell^{\infty}_{\kappa}(\kappa^{++})$ is convex-transitive
for any infinite cardinal $\kappa$.
\end{theorem}
\begin{proof}
Let $\kappa$ be an infinite cardinal and denote $\X=\ell^{\infty}_{\kappa^{+}}(\kappa^{++})\mod \ell^{\infty}_{\kappa}(\kappa^{++})$. 
Fix $x,y\in \S_{\X}$. Pick $u,v\in \ell^{\infty}_{\kappa^{+}}(\kappa^{++})$ such that $x=u+\ell^{\infty}_{\kappa}(\kappa^{++})$ and 
$y=v+\ell^{\infty}_{\kappa}(\kappa^{++})$. 

Thus $||u||\geq 1$ and $||v||\geq 1$. Observe that 
\[|\{\alpha\in \kappa^{++}:\ |u(\alpha)|\geq 1+\frac{1}{k}\}|\leq \kappa\ \mathrm{and}\ |\{\alpha\in \kappa^{++}:\ |v(\alpha)|\geq 1+\frac{1}{k}\}|\leq \kappa\] 
for each $k\in\N$ since $x,y\in \S_{\X}$. Thus, by using the fact that $|\omega_{0}\times \kappa|=\kappa$ we obtain that 
$|\{\alpha\in \kappa^{++}:\ |v(\alpha)|> 1\}|\leq \kappa$ and $|\{\alpha\in \kappa^{++}:\ |u(\alpha)|> 1\}|\leq \kappa$.
Hence $u$ and $v$ above can be chosen such that $u,v\in\S_{\ell^{\infty}_{\kappa^{+}}(\kappa^{++})}$.

Put $A_{i}=\{\alpha\in\kappa^{++}:\ |u(\alpha)|\geq 1-\frac{1}{i}\}$ for $i\in\N$. Observe that
$|A_{i}|=\kappa^{+}$ for $i\in\N$ since $x\in \S_{\X}$. Also note that $|\mathrm{supp}(u)|=|\mathrm{supp}(v)|=\kappa^{+}$.
Fix $B\subset \kappa^{++}\setminus(\mathrm{supp}(u)\cup\mathrm{supp}(v))$ such that $|B|=\kappa^{+}$.
Write $\Gamma=\mathrm{supp}(u)\cup\mathrm{supp}(v)\cup B$.

This finishes the combinatorial part of the argument at hand. The essential point above was to find suitably 
normalized $u,v$ and $\Gamma$. Next we aim to show that when regarding $u,v\in \S_{\ell^{\infty}(\Gamma)}$ it holds that
\begin{equation}\label{eq: vover}
v\in \overline{\conv}(\mathcal{G}_{\ell^{\infty}(\Gamma)}(u))\subset \ell^{\infty}(\Gamma). 
\end{equation}
In such a case one can see similarly as in the proof of Theorem \ref{th prev} that $y\in \overline{\conv}(\mathcal{G}_{\X}(x))$, 
which gives the claim.

To prove \eqref{eq: vover} observe that for each $i\in\N$ there is a bijection $\pi_{i}\colon \Gamma \rightarrow \Gamma$ 
such that $\pi_{i}(A_{i})=\mathrm{supp}(v)$. For each $i\in\N$ define 
$w^{(i)}\colon \Gamma\rightarrow \R$
by 
\[w^{(i)}(\alpha)=\chi_{\mathrm{supp}(v)}(\alpha)\sign(u(\pi_{i}^{-1}(\alpha)))+\chi_{\Gamma\setminus \mathrm{supp}(v)}(\alpha)u(\pi_{i}^{-1}(\alpha)).\]
One can check that 
$v\in \overline{\conv}(\{(\theta(\alpha)w^{(i)}(\alpha))_{\alpha\in \Gamma}|\ \theta\colon \Gamma\rightarrow \{\pm 1\}\ \})$
by applying Fact \ref{Blinfty} on $\mathrm{supp}(v)$ and the (trivial) equality $\frac{w^{(i)}(\alpha)-w^{(i)}(\alpha)}{2}=0$ on
$\Gamma\setminus \mathrm{supp}(v)$ similarly as in \eqref{eq: 4.3*} and \eqref{eq: distz}. 
Finally, Fact \ref{fact2} and the observation that 
$||u(\pi_{i}^{-1})-w^{(i)}||_{\ell^{\infty}(\Gamma)}\leq \frac{1}{i}$ for $i\in\N$ yield \eqref{eq: vover}.
\end{proof}

\section{Final remarks: Constructions for transitive spaces}

Recall that each transitive space contains a separable almost transitive subspace (see \cite[Cor.1.3]{Ca4}, also \cite[Thm. 2.24]{BR2}).
In fact, it turns out below that one can always pass to a suitable transitive subspace of density character at most $2^{\aleph_{0}}$.

\begin{theorem}\label{theorem trans}
Let $\X$ be a transitive space and let $\Y\subset \X$ be a subspace with $\mathrm{dens}(\Y)\leq 2^{\aleph_{0}}$.
Then there exists a closed subspace $\Z\subset\X$ such that
\begin{enumerate}
\item[(1)]{$\Y\subset\Z$}
\item[(2)]{$\mathrm{dens}(\Z)\leq 2^{\aleph_{0}}$}
\item[(3)]{$\Z$ is transitive with respect to subgroup $\mathcal{T}=\{T\in\mathcal{G}_{\X}|\ T(\Z)=\Z\}$ of $\mathcal{G}_{\X}$.}
\end{enumerate}
\end{theorem}

Let us recall the following folklore fact which is proved here for convenience.
\begin{fact}\label{factYd}
Let $(Y,d)$ be a metric space, $\kappa$ a cardinal with $\mathrm{cf}(\kappa)>\omega_{0}$ and 
$\{Y_{\alpha}\}_{\alpha<\kappa}$ an increasing family of closed subsets. Then $\bigcup_{\alpha<\kappa}Y_{\alpha}$ is closed.
\end{fact}
\begin{proof}
Assume to the contrary that there is $x\in \overline{\bigcup_{\alpha<\kappa}Y_{\alpha}}^{d}\setminus \bigcup_{\alpha<\kappa}Y_{\alpha}$.
For each $n<\omega_{0}$ let $\sigma_{n}<\kappa$ be the least ordinal such that $\dist(x,Y_{\sigma_{n}})<\frac{1}{n+1}$.
It follows from the selection of $x$ that $\sup_{n<\omega_{0}}\sigma_{n}=\kappa$, which contradicts $\mathrm{cf}(\kappa)>\omega_{0}$.
\end{proof}

\begin{proof}[Proof of Theorem \ref{theorem trans}]
If $\mathrm{dens}(\X)\leq 2^{\aleph_{0}}$, then the claim holds trivially by putting $\Z=\X$, so let us assume that 
$\mathrm{dens}(\X)>2^{\aleph_{0}}$. We will apply the following fact. If $E$ is a Banach space with $\dim(E)\geq 2$ and 
$\mathrm{dens}(E)\leq 2^{\aleph_{0}}$, then $|\S_{E}|=2^{\aleph_{0}}$. Indeed, clearly $|\S_{E}|\geq 2^{\aleph_{0}}$. 
If $D\subset E$ is dense and $|D|\leq 2^{\aleph_{0}}$, then
of course $D\cap \B_{E}$ is dense in $\B_{E}$ and $|D\cap \B_{E}|\leq 2^{\aleph_{0}}$. Note that for each $e\in \B_{E}$ there is 
a sequence $(d_{n})\subset D\cap \B_{E}$ such that $d_{n}\rightarrow e$ as $n\rightarrow \infty$. Thus, in order to estimate
$|\B_{E}|$, let us estimate the number $|(D\cap \B_{E})^{\omega_{0}}|$ of sequences of $D\cap \B_{E}$ as follows:
\[|(D\cap \B_{E})^{\omega_{0}}|\leq (2^{\aleph_{0}})^{\aleph_{0}}=2^{\aleph_{0}\cdot \aleph_{0}}=2^{\aleph_{0}}.\]
Hence $|\S_{E}|\leq |\B_{E}|\leq 2^{\aleph_{0}}$ and consequently $|\S_{E}|=2^{\aleph_{0}}$. 
Observe that there are 
\begin{equation}\label{eq: SE}
|\S_{E}\times \S_{E}|=|2^{\aleph_{0}}\times 2^{\aleph_{0}}|=2^{\aleph_{0}} 
\end{equation}
many pairs $(x,y)\in \S_{E}\times \S_{E}$.
 
We may assume without loss of generality, possibly by passing to a larger subspace $\Y$, that $\mathrm{dens}(\Y)=2^{\aleph_{0}}$. 
Next we will apply the transitivity of $\X$ in a transfinite induction. 
We will construct increasing sequences $\{\mathcal{T}_{\alpha}\}_{\alpha<2^{\aleph_{0}}}$ and $\{Z_{\alpha}\}_{\alpha<2^{\aleph_{0}}}$, 
which consist of subgroups of $\mathcal{G}_{\X}$ and closed subspaces of $\X$, respectively. 
Put $Z_{0}=\Y$ and $\mathcal{T}_{0}=\{\I_{\X}\}$. For each $0<\alpha<2^{\aleph_{0}}$ we will
construct a subgroup $\mathcal{T}_{\alpha}\subset \mathcal{G}_{\X}$ and a subspace $Z_{\alpha}\subset \X$, 
which satisfy the following conditions:
\begin{enumerate}
\item[(a)]{$|\mathcal{T}_{\alpha}|\leq 2^{\aleph_{0}}$.}
\item[(b)]{$\mathcal{T}_{\beta}\subset \mathcal{T}_{\alpha}$ for $\beta<\alpha$.}
\item[(c)]{For each $x,y\in \S_{\X}\cap \overline{\bigcup_{\beta<\alpha}Z_{\beta}}$ there is $T\in \mathcal{T}_{\alpha}$
such that $T(x)=y$.}
\item[(d)]{$Z_{\alpha}=\overline{\span}(\bigcup \{T(Z_{\beta})|\ T\in\mathcal{T}_{\alpha},\ \beta<\alpha\})$.}
\item[(e)]{$\mathrm{dens}(Z_{\alpha})=2^{\aleph_{0}}$.}
\end{enumerate}
In the case $\alpha=0$ only (a) and (e) apply. 
Suppose that we have obtained this kind of construction for all $\alpha<\gamma$, where $0<\gamma<2^{\aleph_{0}}$. Observe that 
\[\dens\left(\overline{\bigcup_{\alpha<\gamma}Z_{\gamma}}\right)\leq |\gamma\times 2^{\aleph_{0}}|=2^{\aleph_{0}}.\] 
Thus, by \eqref{eq: SE} and the transitivity of $\X$ there is a subset $\mathcal{F}\subset \mathcal{G}_{\X}$
such that for each pair $x,y\in \S_{\X}\cap \overline{\bigcup_{\alpha<\gamma}Z_{\gamma}}$
there is $T\in \mathcal{F}$ such that $T(x)=y$ and $|\mathcal{F}|=2^{\aleph_{0}}$. 
Let $\mathcal{T}_{\gamma}\subset\mathcal{G}_{\X}$ be the subgroup generated by 
the set $\bigcup_{\alpha<\gamma}\mathcal{T}_{\alpha}\cup \mathcal{F}$. Then $\mathcal{T}_{\gamma}$ satisfies conditions 
(b) and (c) (for $\alpha=\gamma$). Note that 
$|\mathcal{T}_{\gamma}|\leq |(\gamma+1)\times 2^{\aleph_{0}}\times \aleph_{0}|=2^{\aleph_{0}}$,
so that $\mathcal{T}_{\gamma}$ satisfies condition (a). Observe that 
$\mathrm{dens}(Z_{\gamma})\leq |\mathcal{T}_{\gamma}\times\gamma\times 2^{\aleph_{0}}\times \aleph_{0}|=2^{\aleph_{0}}$,
so that condition (e) holds, as $\Y\subset Z_{\gamma}$. This completes the induction.

Put $\Z=\bigcup_{\alpha<2^{\aleph_{0}}}\Z_{\alpha}$ and observe that $\dens(\Z)=2^{\aleph_{0}}$.
Note that $\Y\subset\Z$ by the construction. Recall that according to K\"{o}nig's Lemma $\mathrm{cf}(2^{\aleph_{0}})>\aleph_{0}$ 
(see e.g. \cite[Lm. 10.40]{En}). Hence by Fact \ref{factYd} we obtain that
\begin{equation}\label{eq: union}
\Z=\overline{\bigcup_{\alpha<2^{\aleph_{0}}}\Z_{\alpha}}.
\end{equation}

Fix $x,y\in \S_{Z}$ and $T\in \bigcup_{\alpha<2^{\aleph_{0}}}\mathcal{T}_{\alpha}$.
Let $\sigma<2^{\aleph_{0}}$ be such that $T\in \mathcal{T}_{\sigma}$. 
According to \eqref{eq: union}, there is $\beta<2^{\aleph_{0}}$ such that $x,y\in Z_{\beta}$. 
Then $T(x)\in Z_{\max(\sigma,\beta)+1}\subset Z$ by the construction of $Z$. Since $T$ and $x\in \S_{Z}$ were arbitrary,
we get that in fact both $T(Z)\subset Z$ and $T^{-1}(Z)\subset Z$ hold. This means that $T|_{Z}$ is a bijection $Z\rightarrow Z$.

Observe that there is $S\in \mathcal{T}_{\beta+1}$ such that $S(x)=y$ by the construction of $Z$. 
Since $x$ and $y$ were arbitrary and $\bigcup_{\alpha<2^{\aleph_{0}}}\mathcal{T}_{\alpha}\subset \mathcal{T}$, 
we conclude that $Z$ is transitive with respect to $\mathcal{T}$.
\end{proof}
\begin{remark}
Observe that in Theorem \ref{theorem trans} one has some control of the density character of $\Z$, which, even though being a 
superspace of $\Y$, is a still subspace of $\X$. To put Theorem \ref{theorem trans} into perspective, 
recall that every Banach space can be isometrically regarded as 
a subspace of a suitable transitive Banach space (see \cite[Cor. 2.21]{BR2}). Recently Hirvonen and Hyttinen observed
using model theoretic tools \cite[Thm.2.10]{HHpreprint} (see \cite[Thm. 1.13]{HLic} for details) 
that there exist very homogeneous Banach spaces containing all the Banach spaces up to a given density character:
Let $\mu$ be an infinite cardinal. Then there exists a Banach space $\X$ such that
\begin{enumerate}
\item[(1)]{If $\Y$ is any Banach space such that $\dens(\Y)< \mu$, then $\X$ contains an isometric copy of $\Y$.}
\item[(2)]{If closed subspaces $\Y_{1},\Y_{2}\subset\X$ with $\dens(\Y_{1})=\dens(\Y_{2})<\mu$ are mutually isometric, 
then there is $T\in\mathcal{G}_{\X}$ such that $T(\Y_{1})=\Y_{2}$.}
\end{enumerate}

The construction of $\X$ above is not economical in the sense that typically $\mathrm{dens}(\X)$ is much larger than $\mu$. 
Since the $1$-dimensional subspaces $[x]$ and $[y]$ are isometric for all $x,y\in\S_{\X}$, it follows that in (2) 
either $T(x)=y$ or $-T(x)=y$, so that $\X$ is in particular transitive. 
For recent results regarding somewhat economically produced homogeneous spaces, see \cite{Kubis}. 
\end{remark}

\subsection*{Acknowledgements}
This article is part of my Ph.D. research supervised by H.-O. Tylli to whom I am grateful for his careful comments.
This work has been financially supported by the Finnish Graduate School in Mathematical Analysis and Its Applications.

\end{document}